\documentclass[12pt,reqno]{article}

\usepackage[usenames]{color}
\usepackage{amssymb}
\usepackage{graphicx}
\usepackage{amscd}
\usepackage{booktabs}

\usepackage[colorlinks=true,
linkcolor=webgreen,
filecolor=webbrown,
citecolor=webgreen]{hyperref}

\definecolor{webgreen}{rgb}{0,.5,0}
\definecolor{webbrown}{rgb}{.6,0,0}

\usepackage{color}

\usepackage{graphics,amsmath,amssymb}
\usepackage{amsthm}
\usepackage{amsfonts}
\usepackage{latexsym}
\usepackage{epsf}

\newcommand{\seqnum}[1]{\href{http://oeis.org/#1}{\underline{#1}}}

\begin{document}

\theoremstyle{plain}
\newtheorem{theorem}{Theorem}
\newtheorem{corollary}[theorem]{Corollary}
\newtheorem{lemma}[theorem]{Lemma}
\newtheorem{proposition}[theorem]{Proposition}

\theoremstyle{definition}
\newtheorem{definition}[theorem]{Definition}
\newtheorem{question}[theorem]{Question}
\newtheorem{example}[theorem]{Example}
\newtheorem{conjecture}[theorem]{Conjecture}

\theoremstyle{remark}
\newtheorem{remark}[theorem]{Remark}

\begin{center}
\vskip 1cm{\LARGE\bf About Some Relatives of the\\ Taxicab Number}
\vskip 1cm
\large
Viorel Ni\c tic\u a\\
Department of Mathematics\\
West Chester University of Pennsylvania\\
West Chester, PA 19383\\
USA\\
\href{mailto:vnitica@wcupa.edu}{\tt vnitica@wcupa.edu} \\
\end{center}

\vskip .2 in

\begin{abstract} The taxicab number, $1729$, is the smallest number that can be written as a sum of two cubes in two different ways. It also has the following property: if we add its digits we obtain $19$. The number obtained from $19$ reversing the order of its digits is $91$. If we multiply $19$ by $91$ we obtain again $1729$. In the paper we study various generalizations of this property.
\end{abstract}

\section{Introduction}\label{sec:1}

The \emph{taxicab number,} $1729$, became well known due to a discussion between Hardy and Ramanujan \cite{Hardy}. It is the smallest positive integer that can be written in two ways as a sum of two cubes: $1^3+12^3$ and $9^3+10^3$. The number $1729$ also has a less well known property: if we add its digits we obtain $19$; multiplying $19$ by $91$, the number obtained from 19 reversing the order of its digits, we obtain again $1729$. It is not hard to show that the set of integers with this property is finite and equal to $\{1, 81, 1458, 1729\}$.

In a conversation that the author had with his colleague, Professor Shiv Gupta, Shiv asked if the second property can be generalized. One replaces the sum of the digits of an integer by the sum of the digits times an integer multiplier and then multiplies the product by the number obtained reversing the order of the digits in the product. The taxicab number becomes a particular example with multiplier 1. A computer search produced a large number of examples with larger multiplier. There are 23 integers less than 10000 having this property; see sequence  \seqnum{A305131} in the OEIS \cite{online}. For example, 2268 has multiplier 2. The sum of the digits is $18$, one has $18\times 2=36$, and $36\times 63=2268$.

One may replace the last product in the above procedure by a sum. A computer search showed that there are numbers that have the property for sums. There are 264 integers less than 10000 having the property; see sequence \seqnum{A305130}  in the OEIS \cite{online}. For example, $121212$ has multiplier $6734$. The sum of the digits is $9$, one has $9\times 6732=60606$, and $60606+60606=121212$.

The paper is dedicated to the study of these properties. After the paper was submitted for publication we learned from the editor that our work may be related to the study of Niven (or Harshad) numbers. These are numbers divisible by the sum of their decimal digits. Niven numbers have been extensively studied as one can see for instance from Cai \cite{C}, Cooper and Kennedy \cite{CK}, De Koninck and Doyon \cite{KD}, Grundman \cite{G}. One of the classes of integers we study, that of multiplicative Ramanujan-Hardy numbers, is a subclass of the class of Niven numbers. Of interest are also $q$-Niven numbers, which are numbers divisible by the sum of their base $q$ digits. See, for example, Fredricksen, Iona\c scu, Luca, and St\u anic\u a \cite{FILS}. Some other variants of Niven numbers can be found in Boscaro \cite{B1} and Bloem \cite{B2}.

\section{Statements of the main results}\label{sec:1-bis}

In what follows let $b\ge 2$ be an arbitrary numeration base.

\begin{definition} If $N$ is a positive integer written in base $b$, we call \emph{reversal} of $N$ and let $N^R$ denote the integer obtained from $N$ by writing its digits in reverse order.
\end{definition}

We observe that addition and multiplication are independent of the numeration base. The operation of taking the reversal is not.

Let $s_b(N)$ denote the sum, done in base 10, of the base $b$ digits of an integer $N$.

\begin{definition}\label{def:1} A positive integer $N$ written in base $b$ is called \emph{$b$-additive Ramanujan-Hardy number,} or simply $b$-ARH number, if there exists a positive integer $M$, called \emph{additive multiplier}, such that
\begin{equation}\label{eq:1}
N=Ms_b(N)+(Ms_b(N))^R,
\end{equation}
where $(Ms_b(N))^R$ is the reversal of base $b$-representation of $Ms_b(N)$.
\end{definition}

\begin{definition}\label{def:2} A positive integer $N$ written in base $b$ is called \emph{$b$-multiplicative Ramanujan-Hardy number,} or simply $b$-MRH number, if there exists a positive integer $M$, called \emph{multiplicative multiplier}, such that
\begin{equation}\label{eq:2}
N=Ms_b(N)\cdot (Ms_b(n))^R,
\end{equation}
where $(Ms_b(N))^R$ is the reversal of base $b$-representation of $Ms_b(N)$.
\end{definition}

To simplify the notation, let $s(N)$, ARH, MRH denote $s_{10}(N)$, 10-ARH, 10-MRH.

While  $b$-MRH numbers are $b$-Niven numbers, $b$-Niven numbers are not necessarily $b$-MRH numbers.

\begin{example} The number $[144]_7$ is a $7$-Niven number but not a $7$-MRH number.
\end{example}

We observe that the notions of $b$-ARH and $b$-MRH numbers are dependent on the base.

\begin{example} The number $[12]_{10}$ is an ARH number, but $[12]_{9}$ is not a $9$-ARH number. The number $[81]_{10}$ is an MRH number, but $[81]_{9}$  is not a $9$-MRH number.
\end{example}

Once these notions are introduced, several natural questions arise.

\begin{question}\label{q:1} Do there exist infinitely many $b$-ARH numbers?
\end{question}

\begin{question}\label{q:2} Do there exist infinitely many $b$-MRH numbers?
\end{question}

\begin{question}\label{q:3} Do there exist infinitely many additive multipliers?
\end{question}

\begin{question}\label{q:4} Do there exist infinitely many multiplicative multipliers?
\end{question}

In what follows, if $x$ is a string of digits, we let $(x)^{\land k}$ denote the base 10 integer obtained by repeating $x$ $k$-times. We also let $[x]_b$ denote the value of the string $x$ in base $b$.

The following example gives an explicit positive answer to Question \ref{q:1} if $b=10$.

\begin{example}\label{thm:1} Consider the numbers
\begin{equation}\label{eq:3}
N_k=(12)^{\land 3^k},
\end{equation}
where $k$ is a positive integer. All numbers $N_k$ are ARH numbers and Niven numbers.
In particular,  there exist infinitely many Niven numbers with no digit equal to zero.
\end{example}

\begin{example}\label{r:1} If we allow zero digits an infinity of $b$-MRH numbers is given by $\{[1(0)^{\land k}]_b\vert k\in \mathbb{N}\}$. Last example has the unpleasant feature that the apparent multiplicative multiplier of each $b$-MRH numbers is the number itself and the search for other multipliers is dependent on the base. In order to avoid trivial considerations, we consider from now on only examples of $b$-ARH and $b$-MRH numbers that have many digits different from zero.
\end{example}

It follows from the proof of Example \ref{thm:1} that $Ms(N_k)=(Ms(N_k))^R$. The following theorem gives an example in which
it is clear from the proof that $Ms_b(N_k)\not =(Ms_b(N_k))^R$ for an arbitrary even base $b$. One can read from the proof the explicit base $b$ expansion of the multipliers. Counting the multipliers shows that the set of multipliers of a $b$-ARH number $N$ can grow exponentially in terms of the number of digits of $N$.

\begin{theorem}\label{thm:11} Consider the numbers
\begin{equation}\label{eq:3**}
N_k=[(1)^{\land k}]_b,
\end{equation}
where $b$ is even, $k=[1(0)^{\land p}]_b, p\ge 1$, $p$ an arbitrary natural number. All numbers $N_k$ are $b$-ARH numbers and not $b$-Niven numbers.

Each $N_k$ has a subset of additive multipliers of cardinality $2^{\frac{k-2p}{2}}$ consisting of all integers $[(1)^{\land p}I]_b$, where $I$ is a sequence of $0$ and $1$ of length $k-2p$ in which no two digits symmetric about the center of the sequence are identical.
\end{theorem}

\begin{example} We show an example that illustrates the results in Theorem \ref{thm:11}. Assume that $b=2$, $k=16=[10000]_2$, and $p=4$. Then $N_{16}=[(1)^{\land 16}]_2$ and $s_2(N_{16})=2^4=[10000]_2$. The following $16=2^{\frac{16-2\cdot 4}{2}}$ numbers are additive multipliers of $N_{16}$:
\begin{equation*}
\begin{gathered}
\ [111100001111]_2,\ [111100010111]_2,\ [111100101011]_2,\ [111100111100]_2,\\
\ [111101001101]_2,\ [111101010101]_2,\ [111101101001]_2,\ [111101110001]_2,\\
\ [111110001110]_2,\ [111110010110]_2,\ [111110101010]_2,\ [111110110010]_2,\\
\ [111111001100]_2,\ [111111010100]_2,\ [111111101000]_2,\ [111111110000]_2.
\end{gathered}
\end{equation*}
\end{example}

\begin{remark} The numbers $N_k$ may have other multipliers, besides those listed in Theorem \ref{thm:11}. The growth of the set of multipliers can be larger than that shown in Theorem \ref{thm:11} and depends on the numeration base; see Theorem \ref{thm:larger-growth}. Nevertheless, for $b=2$ there are no other multipliers of $N_k$ besides those listed in Theorem \ref{thm:11}. We observe that the numbers $N_k$ from Theorem \ref{thm:11} have an even number of digits and the numbers $N_k$ from Theorem \ref{thm:larger-growth} have an odd number of digits.
\end{remark}

\begin{theorem}\label{thm:larger-growth} Consider the numbers
\begin{equation}\label{eq:3larger-growth}
N_k=[(1)^{\land p}(10)^{\land k-2p}0(1)^{\land p}]_b,
\end{equation}
where $b$ is even and $k=[1(0)^{\land p}]_b, p\ge 1,$ $p$ arbitrary natural number. All numbers $N_k$ are $b$-ARH numbers and not $b$-Niven numbers.

For each $N_k$ the set of additive multipliers has cardinality $(b-1)^\frac{k-2p}{2}$ and consists of all integers $[(1)^{\land p}I0]_b$, where $I$ is a concatenation of $k-2p$ two digits strings of type $0\alpha, \alpha\not = 0$, in which any pair of nonzero digits symmetric about the center of $I0$ have their sum equal to $b$.
\end{theorem}

\begin{example} We show an example that illustrates the results in Theorem \ref{thm:larger-growth}. Assume that $b=4$, $k=4=[10]_4$, and $p=1$. Then $N_{4}=[1101001]_4$ and $s_4(N_{4})=4=[10]_4$. The following $3=3^\frac{4-2\cdot 1}{2}$ numbers are additive multipliers of $N_{4}$:
\begin{equation*} \ [102020]_4, \  [101030]_4, \ [103010]_4.
\end{equation*}
\end{example}

The following corollary of Theorem \ref{thm:2} gives a partial answer to Question \ref{q:3}.

\begin{corollary} If $b$ is even there exist infinitely many additive multipliers. Moreover, there exists infinitely many $b$-ARH numbers that have at least two additive multipliers.
\end{corollary}

The numbers $N_k$ from Theorems \ref{thm:11} and \ref{thm:larger-growth} are not $b$-MRH numbers.

\begin{question} Do there exist infinitely many $b$-MRH numbers that have at least two multiplicative multipliers?
\end{question}

\begin{corollary} If $b$ is even there exist infinitely many $b$-ARH numbers that are not $b$-MRH.
\end{corollary}

Motivated by the results in Theorems \ref{thm:11} and \ref{thm:larger-growth} and by the examples of ARH and MRH numbers shown in Sections \ref{sec:7} and \ref{sec:8}, we introduce the following notions.

\begin{definition} If $N$ is a $b$-ARH number, let the \emph{multiplicity} of $N$ be the cardinality of the corresponding set of additive multipliers.
\end{definition}

\begin{definition} If $N$ is a $b$-MRH number, let the \emph{multiplicity} of $N$ be the cardinality of the corresponding set of multiplicative multipliers.
\end{definition}

Theorem \ref{thm:11} has the following corollary.

\begin{corollary} The multiplicity of $b$-ARH numbers is unbounded for any even base.
\end{corollary}

\begin{question}\label{q:13} Is the multiplicity of $b$-MRH numbers bounded?
\end{question}

\begin{remark}\label{r:2} For Questions \ref{q:2} and \ref{q:4} we do not have an answer with $b$-MRH numbers having all digits different from zero. See Theorem \ref{thm:mrhexample} for an infinity of $b$-MRH numbers with half of the digits different from zero. No prime number can be an MRH number. Note that no integer with  two prime factors in the prime factorization can be an MRH number. Such an MRH number has the multiplier equal to $1$ and among the MRH numbers with multiplier $1$ none has two factors in the prime factorization.
\end{remark}

The following theorem shows an infinity of $b$-Niven numbers that are not $b$-MRH numbers.

\begin{theorem}\label{thm:niven-not-mrh} Let $b\ge 2$ be a numeration base. For $n$ not divisible by $b-1$ define
\begin{equation*}
R_n=\frac{b^n-1}{b-1}=[(1)^{\land n}]_{b}, n\ge 1.
\end{equation*}
Then $(b-1)nR_n$ is a $b$-Niven number that is not  a $b$-MRH number.
\end{theorem}

For $b$-ARH numbers one has the following result.

\begin{theorem}\label{thm:noARH} There exist infinitely many integers that are not $b$-ARH numbers.
\end{theorem}

The following Theorem gives a partial answer to Question \ref{q:2}.

\begin{theorem}\label{thm:mrhexample} Let $b$ odd and $k\ge 2$. Then the numbers
\begin{equation}\label{eq:mrhex}
N_k=[(b-1)^{\land 2^{k-1}-1}(b-2)(0)^{\land 2^{k-1}-1}1]_b
\end{equation}
are $b$-MRH numbers and $s_b(\sqrt{N_k})=s_b(N_k)$.

Moreover, if $b\equiv 3 \pmod {4}$ then $\sqrt{N_k}$ is itself a $b$-Niven number.
\end{theorem}

\begin{example} We illustrate the result in Theorem \ref{thm:mrhexample}.
\begin{itemize}
\item For $b=3, k=2$ we get $N_2=[2101]_3$ which is a $3$-MRH number. Then  $\sqrt{[2101]_3}=[22]_3$, $s_3([2101]_3)=s_3([22]_3)=4$ and $[22]_3$ is a $3$-Niven number.
\item For $b=5, k=2$ we get $N_2=[4301]_5$ which is a $5$-MRH number. Then $\sqrt{[4301]_5}=[44]_5$, $s_5([4301]_5)=s([44]_5)=8$ and $[44]_5$ is a $5$-Niven number.
\item For $b=17, k=5$, $N_5$ is a $17$-MRH number, but $\sqrt{N_5}$ is not a $17$-Niven number.
\item For $b=7, k=2$ we get $N_2=[6501]_7$ which is a $7$-MRH number. Then $\sqrt{[6501]_7}=[66]_7$, $s_7([6501]_7)=s_7([66]_7)=12$ and $[66]_7$ is a $7$-Niven number.
\end{itemize}

Third item shows that the congruence condition in Theorem \ref{thm:mrhexample} is necessary. Second item shows that $\sqrt{N_k}$ may be a $b$-Niven number even without this condition.
\end{example}

The following corollary of Theorem \ref{thm:mrhexample} gives a partial positive answer to Question \ref{q:4}.

\begin{corollary} If $b$ is odd there exist infinitely many multiplicative multipliers.
\end{corollary}

We show two unexpected corollaries of the proof of Theorem \ref{thm:mrhexample}.

\begin{corollary} If $b$ is odd there exist infinitely many $b$-MRH numbers that are perfect squares.
\end{corollary}

\begin{corollary}\label{coro:number} If $b\equiv 3\pmod{4}$ there exists an infinity of $b$-MRH numbers $N$ for which $\sqrt{N}$ is a $b$-Niven number and for which $s_b(N)=s_b(\sqrt{N})$.
\end{corollary}

The following notions of high degree $b$-Niven numbers are motivated by Corollary \ref{coro:number}, which provides plenty of examples.

\begin{definition} An integer $N$ is called \emph{quadratic $b$-Niven number} if $N$ and $N^2$ are $b$-Niven numbers.
If in addition $s_b(N)=s_b(N^2)$ then $N$ is called \emph{strongly quadratic $b$-Niven number}.
\end{definition}

The study of high degree $b$-Niven numbers is continued in Ni\c tic\u a \cite{N}. We show that for each degree there exists an infinity of bases in which $b$-Niven numbers of that degree appear.

We show in Sections \ref{sec:7} that $6$ is not an additive multiplier for base $10$ and ARH numbers without zero digits, and that $9$ is not an additive multiplier for base $10$. We show in Section \ref{sec:8} that $3$ is not a multiplicative multiplier for base $10$. We do not know how to answer the following questions for any base.

\begin{question}\label{q:6} Do there exist infinitely many integers that are not additive multipliers?
\end{question}

\begin{question}\label{q:7} Do there exist infinitely many integers that are not multiplicative multipliers?
\end{question}

In what follows let $\lfloor x \rfloor$ denote the integer part, let $\ln x$ denote the natural logarithm and let $\log_b x$ denote base $b$ logarithm of the positive real number $x$.

The following theorems give bounds for the number of digits in a $b$-ARH number in terms of the multiplier.

\begin{theorem}\label{thm:2} Let $N$ be a $b$-ARH number with $k$ digits and  additive multiplier $M$. Then
\begin{displaymath}
k\leq \begin{cases}
	M+2, & \text{ if } b\ge 4;\\
	M+3, & \text{ if } b=2 \text{ or } b=3.
\end{cases}
\end{displaymath}
\end{theorem}

\begin{corollary} For fixed additive multiplier $M$ and base $b$, the set of $b$-ARH numbers with multiplier $M$ is finite.
\end{corollary}

\begin{theorem}\label{thm:2-strong} Let $N$ be a $b$-ARH number with $k$ digits and additive multiplier $M$. Under any of the following assumptions:
\begin{itemize}
\item $b \ge 10$ and $M\ge b^6;$
\item $3\le b \le 9$ and $M\ge b^7;$
\item $b=2$ and $M\ge b^8,$
\end{itemize}
one has
\begin{equation}\label{eq:4-Strong}
\begin{gathered}
k\le 2\lfloor \log_b M \rfloor.
\end{gathered}
\end{equation}
\end{theorem}

The following theorems give bounds for the number of digits in a $b$-MRH number in terms of the multiplier.

\begin{theorem}\label{thm:3} Let $N$ be a $b$-MRH number with $k$ digits and multiplicative multiplier $M$. Then
\begin{displaymath}
k\leq \begin{cases}
	M+4, & \text{ if } b\ge 6;\\
	M+5, & \text{ if } b= 5;\\
     M+7, & \text{if } 2\le b\le 4.
\end{cases}
\end{displaymath}
\end{theorem}

Theorem \ref{thm:3} shows that a MRH number with multiplicity $1$ can have at most 5 digits. A computer search shows that the set of all such numbers is indeed $\{1,81, 1458, 1729\}$.

\begin{corollary} For fixed multiplicative multiplier $M$ and base $b$, the set of $b$-MRH numbers with multiplier $M$ is finite.
\end{corollary}

\begin{theorem}\label{thm:3-new} Let $N$ be a $b$-MRH number with $k$ digits and multiplicative multiplier $M$. Under any of the following assumptions:
\begin{itemize}
\item $b \ge 9$ and $M\ge b^9;$
\item $5\le b \le 8$ and $M\ge b^{10};$
\item $b=4$ and $M\ge b^{11};$
\item $b=3$ and $M\ge b^{12};$
\item $b=2$ and $M\ge b^{16};$
\end{itemize}
one has
\begin{equation}\label{eq:5-new}
k\le 3\lfloor \log_b M\rfloor.
\end{equation}
\end{theorem}

We summarize the rest of the paper. Example \ref{thm:1} is proved in Section \ref{sec:2}, Theorem \ref{thm:11} is proved in Section \ref{sec:2bis}, Theorem \ref{eq:3larger-growth} is proved in Section \ref{sec:3larger-growth}, Theorem \ref{thm:noARH} is proved in Section \ref{sec:noARH}, Theorem \ref{thm:niven-not-mrh} is proved in Section \ref{sec:niven-not-mrh}, Theorem \ref{thm:mrhexample} is proved in Section \ref{sec:mrhexample}, Theorem \ref{thm:2} is proved in Section \ref{sec:3}, Theorem \ref{thm:2-strong} is proved in Section \ref{sec:4}, Theorem \ref{thm:3} is proved in Section \ref{sec:5}, and Theorem \ref{thm:3-new} is proved in Section \ref{sec:6}. In Section \ref{sec:7} we show examples of ARH numbers and ask additional questions and in Section \ref{sec:8} we show examples of MRH numbers and ask additional questions. In Section \ref{sec:9} we describe an approach to Question \ref{q:2} if $b=10$.

\section{Proof of Example \ref{thm:1}}\label{sec:2}

\begin{proof} One obtains a formula for $N_k$ by adding two geometric series.
\begin{equation}\label{eq:5}
\begin{gathered}
N_k=10^{2\cdot 3^k-1}+10^{2\cdot 3^k-3}+\ldots +10\\
+2(10^{2\cdot 3^k-2}+10^{2\cdot 3^k-4}+\ldots +1)\\
=12\cdot \frac{10^{2\cdot 3^k}-1}{99}=4\cdot \frac{10^{2\cdot 3^k}-1}{33}.
\end{gathered}
\end{equation}

Note that $s(N_k)=3^{k+1}$. We show by induction that $s(N_k)$ divides $N_k$. The case $k=0$ gives $s(N_0)=3$ which divides $N_0=12$. Assume that for fixed $k$, $s(N_k)$ divides $N_k$.

\begin{equation}\label{eq:6}
\begin{gathered}
N_{k+1}=4\cdot \frac{10^{2\cdot 3^{k+1}}-1}{33}=4\cdot \frac{(10^{2\cdot 3^k})^3-1^3}{33}\\
=4\cdot \frac{10^{2\cdot 3^k}-1}{33}\cdot (10^{4\cdot 3^k}+10^{2\cdot 3^k}+1)=N_k\cdot (10^{4\cdot 3^k}+10^{2\cdot 3^k}+1),
\end{gathered}
\end{equation}
which is clearly divisible by $s(N_{k+1})=3^{k+2}$ due to $N_k$ divisible by $s(N_k)=3^{k+1}$ and $10^{4\cdot 3^k}+10^{2\cdot 3^k}+1$ divisible by $3$. Therefore $s(N_k)$ divides $N_k$ and $N_k$ is a Niven number.

Observe now that $N_k/2=(N_k/2)^R$. It follows from \eqref{eq:3} and the fact that $N_k$ is divisible by $s(N_k)=3^{k+1}$ that $N_k/2$ is divisible by $s(N_k)$. We conclude that $N_k$ is an ARH number with additive multiplier $M=N_k/(2s(N_k))$.
\end{proof}

\section{Proof of Theorem \ref{thm:11}}\label{sec:2bis}
\begin{proof} Let $N_k=[(1)^{\land k}]_b$ where $k$ is even and $k=[1(0)^{\land p}]_b, p\ge 1$, $p$ arbitrary natural number. Then $s_b(N_k)=[1(0)^{\land p}]_b$. Let $M=[(1)^{\land p}I]_b$, where $I$ is a string of $0$ and $1$ of length $k-2p$ in which no two digits symmetric about the center of the sequence are identical. Note that $M^R=[(I)^R(1)^{\land p}]_b$. The following calculation shows that $N_k$ is a $b$-ARH number. Note that $I+(I)^R=[(1)^{\land k-2p}]_b$.

\begin{equation*}
\begin{gathered}
s_b(N_k)\cdot M+(s_b(N_k)\cdot M)^R\\
=[1(0)^{\land p}]_b\cdot [(1)^{\land p}I]_b+([1(0)^{\land p}]_b\cdot [(1)^{\land p}I]_b)^R\\
=[(1)^{\land p}I(0)^{\land p}]_b+([(1)^{\land p}I(0)^{\land p}]_b)^R\\
=[(1)^{\land p}I(0)^{\land p}]_b+[(0)^{\land p}(I)^R(1)^{\land p}]_b=[(1)^{\land k}]_b=N_k.
\end{gathered}
\end{equation*}

In order to count the multipliers, observe that the length of the string $I$ is $k-2p$. If we know half of its digits we can find the other half using that no two digits symmetric about the center of the string are identical. The number of strings of $0$ and $1$ of length $\frac{k-2p}{2}$ is $2^{\frac{k-2p}{2}}$. Finally, observe that $N_k$ is not divisible by $s_b(N_k)$, so $N_k$ is not a $b$-Niven number..
\end{proof}

\section{Proof of Theorem \ref{thm:larger-growth}}\label{sec:3larger-growth}

\begin{proof} Let $N_k=[(1)^{\land p}(10)^{\land k-2p}0(1)^{\land p}]_b$ where $b$ is even and $k=[1(0)^{\land p}]_b, p\ge 1$. Then $s_b(N_k)=[1(0)^{\land p}]_b$. Let $M=[(1)^{\land p}I0]_b$. Note that $M^R=[0(I)^R(1)^{\land p}]_b$. The following calculation shows that $N_k$ is a $b$-ARH number. Note that $I0+0(I)^R=[(10)^{\land k-2p}0]_b$.

\begin{equation*}
\begin{gathered}
s_b(N_k)\cdot M+(s_b(N_k)\cdot M)^R\\
=[1(0)^{\land p}]_b\cdot [(1)^{\land p}I0]_b+([1(0)^{\land p}]_b\cdot [(1)^{\land p}I0]_b)^R\\
=[(1)^{\land p}I0(0)^{\land p}]_b+([(1)^{\land p}I0(0)^{\land p}]_b)^R\\
=[(1)^{\land p}I0(0)^{\land p}]_b+[(0)^{\land p}0(I)^R(1)^{\land p}]_b=[(1)^{\land p}(10)^{\land k-2p}0(1)^{\land p}]_b=N_k.
\end{gathered}
\end{equation*}

In order to count the multipliers, observe that the number of nonzero digits in the string $I0$ is $k-2p$. If we know half of the nonzero digits we can find the other half using that no two digits symmetric about the center of the string $I0$ are identical. There are $\frac{k-2p}{2}$ positions to be filled and each one can be filled in $b-1$ ways. To show that there are no other multiplier it is enough to prove, using induction on length, that the string $[(10)^{\land k-2p}0]_b$ cannot be written as a sum of a string $J$ and its reversal except if $J=I0$, where $I$ is as above.
Finally, observe that $N_k$ is not divisible by $s_b(N_k)$, so $N_k$ is not a $b$-Niven number.
\end{proof}

\section{Proof of Theorem \ref{thm:niven-not-mrh}}\label{sec:niven-not-mrh}

\begin{proof} McDaniels proved \cite[Theorem 2]{D} that if $b=10$ and $m \le 9R_n$ then $s(9mR_n)=9n$. The proof is valid in any base $b$ and follows readily upon writing $m$ as:
\begin{equation}
m=\sum_{i=0}^k a_ib^i, k<n.
\end{equation}
It gives that if $m\le (b-1)R_n$ then $s_b((b-1)mR_n)=(b-1)n$.
If $m=n$ one has  $s_b((b-1)nR_n)=(b-1)n$, which shows that $(b-1)nR_n$ is a $b$-Niven number.
By contradiction, assume that $(b-1)nR_n$ is a $b$-MRH number with multiplier $M$. It follows that:
\begin{equation}\label{eq:contrad}
(b-1)nM((b-1)nM)^R=(b-1)nR_n.
\end{equation}

We recall that a base $b$ number is divisible by $b-1$ if the sum of its base $b$ digits is divisible by $b-1$. The divisibility test and $b-1\not \vert n$ imply that $b-1\not\vert R_n$, but $b-1\vert ((b-1)nM)^R$. As $b-1\not\vert n$, there are at least two factors of $b-1$ in the factorization of the left hand side of \eqref{eq:contrad} and only one factor of $b-1$ in the right hand side of \eqref{eq:contrad}. This gives a contradiction.
\end{proof}

\section{Proof of Theorem \ref{thm:noARH}}\label{sec:noARH}

\begin{proof} A $b$-ARH number is a sum of an integer and its reversal. In order to prove the theorem it is enough to show that there exist infinitely many integers that are not a sum of an integer and its reversal. There are $b^{k}-b^{k-1}=b^{k-1}(b-1)$ base $b$ $k$-digit numbers. Those of type $N+N^R$, either have $N=[a_ka_{k-1}\cdots a_2a_1]_b$ with $a_k+a_1\le b-1$, or have $N$ with $k-1$ digits. There are $\frac{b(b-1)}{2}\cdot b^{k-2}$ $k$-digit numbers with $a_k+a_1\le b-1$ and there are $b^{k-1}-b^{k-2}$ $(k-1)$-digit numbers. Overall, we have
\begin{equation*}
\frac{b(b-1)}{2}\cdot b^{k-2}+(b^{k-1}-b^{k-2})=b^{k-1}\left ( \frac{b+1}{2}\right )-b^{k-2}
\end{equation*}
$k$-digit numbers of type $N+N^R$. Hence there are
\begin{equation}\label{eq:new-add1}
b^k-b^{k-1}-\left ( b^{k-1}\left ( \frac{b+1}{2}\right )-b^{k-2}\right )=b^{k-1}\left ( \frac{b-3}{2}\right )+b^{k-2}
\end{equation}
$k$-digit numbers that are not of type $N+N^R$. The right hand side of equation \eqref{eq:new-add1} has limit $\infty$ as $k\to \infty$ for $b\ge3$ and this ends the proof if $b\ge 3$. If $b=2$, consider the numbers $[(1)^{\land k}]_2$. These are not ARH-numbers if $k$ is odd.
\end{proof}

\section{Proof of Theorem \ref{thm:mrhexample}}\label{sec:mrhexample}

\begin{proof} As $\gcd(b,2)=1$ Euler's Theorem implies that $2^k$ divides $b^{\phi(2^k)}-1$. Clearly $b-1$ also divides $b^{\phi(2^k)}-1$. Assume that $\text{gcd}(2^k,b-1)=2^\ell$. Then $2^{k-\ell}(b-1)$ divides $b^{\phi(2^k)}-1=b^{2^{k-1}}-1$. Consider the product
\begin{equation*}
(b^{2^{k-1}}-1)^2=b^{2\cdot2^{k-1}}-2b^{2^{k-1}}+1.
\end{equation*}
The product is divisible by $2^{k-1}(b-1)$, written in base $b$ equals $N_k$, and $s_b(N_k)=2^{k-1}(b-1)$. We conclude that $N_k$ is a $b$-MRH number.

To finish the proof of the theorem observe that if $b\equiv 3\pmod{4}$ then $\text{gcd}(2^k,b-1)=2$. Therefore $2^{k-1}(b-1)$ divides $b^{2^k}-1=[(b-1)^{2^k-1}]_b$. Finally \begin{equation*}
s_b(\sqrt{N_k})=s_b([(b-1)^{2^k-1}]_b)=2^{k-1}(b-1)\vert b^{2^k}-1=\sqrt{N_k}.
\end{equation*}
\end{proof}

\section{Proof of Theorem \ref{thm:2}}\label{sec:3}

\begin{proof} As $N$ has $k$ digits one has that:
\begin{equation}\label{eq:7}
N\ge b^{k-1}.
\end{equation}

The largest possible value for $s_b(N)$ is $(b-1)k$. We observe that reversing the order of the digits in an integer increases its value by at most $b$ times. One has that:
\begin{equation}\label{eq:8}
Ms_b(N)+(Ms_b(N))^R\le (b^2-1)kM.
\end{equation}

Combining equations \eqref{eq:1}, \eqref{eq:7}, \eqref{eq:8} one has that:
\begin{equation}\label{eq:9}
b^{k-1} \le (b^2-1)kM.
\end{equation}

We prove by induction on the variable $k$ that:
\begin{equation}\label{eq:10}
b^{k-1} > (b^2-1)kM, \text{ for }k\ge M+3, M\ge 1, b\ge 4,
\end{equation}
which combined with \eqref{eq:9} gives a contradiction and ends the proof of Theorem \ref{thm:2} for $b\ge 4$.

In the first step $k=M+3$. The statement in \eqref{eq:10} becomes
\begin{equation}\label{eq:11}
b^{M+2} > (b^2-1)(M^2+3M), \text{ for }M\ge 1, b\ge 4.
\end{equation}

We prove \eqref{eq:11} by induction on the variable $M$. In the initial step $M=1$ and one has
\begin{equation*}
b^3 > 4(b^2-1) \Leftrightarrow b^2(b-4)+4 > 0,
\end{equation*}
which is clearly true for $b\ge 4$.

Now assume that \eqref{eq:11} is true for $M$ and prove it for $M+1$. Using the induction hypothesis one has that:
\begin{equation}\label{eq:12}
b^{M+3}=b\cdot b^{M+2}> b\cdot (b^2-1)(M^2+3M).
\end{equation}

In order to finish the proof by induction, we still need to check that:
\begin{equation}\label{eq:13}
b\cdot (b^2-1)(M^2+3M)\ge (b^2-1)\left ( (M+1)^2+3(M+1)\right ).
\end{equation}

After simplifications, \eqref{eq:13} becomes
\begin{equation}\label{eq:14}
(b-1)M^2+(3b-5)M-4\ge 0.
\end{equation}

As the left hand side of \eqref{eq:14} is larger than $M^2+4M-4$, which is positive if $M\ge 2$, we conclude that \eqref{eq:14} is true for all $M\ge 1$ and finish the proof of \eqref{eq:11}.

We continue with the general step in the proof of \eqref{eq:10}. By induction:
\begin{equation}\label{eqn:ind1}
b^k=b\cdot b^{k-1}>b(b^2-1)kM.
\end{equation}

We still need to check that
\begin{equation}\label{eqn:ind2}
b(b^2-1)kM\ge (b^2-1)k(M+1),
\end{equation}
which is obvious and finishes the proof of \eqref{eq:11} and that of Theorem \ref{thm:2} for base $b\ge 4$.

Now assume $b=3$. Equation \eqref{eq:9} is still valid.

We prove by induction on the variable $k$ that:
\begin{equation}\label{eq:some-new-1}
b^{k-1} > (b^2-1)kM, \text{ for }k\ge M+4, M\ge 1.
\end{equation}
Equations \eqref{eq:9} and \eqref{eq:some-new-1} give a contradiction that finishes the proof of the theorem.

If $k=M+4$ one has that:
\begin{equation}\label{eq:some-new-2}
b^{M+3} > (b^2-1)(M^2+4M), \text{ for } M\ge 1,
\end{equation}
which we prove by induction on $M$.

The case $M=1$ is true. We assume \eqref{eq:some-new-2} true for $M$ and prove it for $M+1$. By induction one has that:
\begin{equation*}
b^{M+4}=b\cdot b^{M+3} > b(b^2-1)(M^2+4M).
\end{equation*}

To finish the proof of \eqref{eq:some-new-2} we still need to check that:
\begin{equation}\label{eq:1corec}
b(b^2-1)(M^2+4M)\ge (b^2-1)\left ( (M+1)^2+4(M+1)\right ),
\end{equation}
which simplifies to $(b-1)M^2+(4b-6)M-5\ge 0$ and is true for $M\ge 1, b=3$.

The rest of the proof of \eqref{eq:some-new-1} follows from \eqref{eqn:ind1} and \eqref{eqn:ind2}.

Assume $b=2$. Equation \eqref{eq:9} is still valid.

We prove by induction on the variable $k$ that:
\begin{equation}\label{eq:some-new-6}
b^{k-1} > (b^2-1)kM, \text{ for }k\ge M+4, M\ge 3.
\end{equation}
Equations \eqref{eq:9} and \eqref{eq:some-new-6} give a contradiction that ends the proof of the theorem for $b=2, M\ge 3$.

If $k=M+4$ one has that:
\begin{equation}\label{eq:some-new-60}
2^{M+3} > 3(M^2+4M), \text{ for } M\ge 3,
\end{equation}
which we prove by induction on $M$.

The case $M=3$ is true. Assume now that \eqref{eq:some-new-60} is true for $M$ and prove it for $M+1$.

\begin{equation*}
2^{M+4}=2\cdot 2^{M+3} > 6(M^2+4M).
\end{equation*}

To finish the proof we still need to check that:
\begin{equation*}
6(M^2+4M)\ge 3\big ( (M+1)^2+4M\big ).
\end{equation*}
The equation simplifies to $M^2+2M-1\ge 0$ and it is true for $M\ge 3$.

To finish the proof of the theorem if $b=2$, it remains to discuss the cases $M=1, M=2$.

Let $M=1$. If $k\le 4$ the theorem is trivially true, so assume $k\ge 5$.
Let $N$ be a $2$-ARH number with $k$ digits and $M=1$. Then $s_2(N)\le k$ and $N\ge 2^{k-1}$. This implies
\begin{equation}
2^{k-1}\le 3k.
\end{equation}
One shows that $3k<2^{k-1}$ for $k\ge 5$ and gets a contradiction.

Let $M=2$. If $k\le 4$ the theorem is trivially true, so assume $k\ge 5$.
Let $N$ be a $2$-ARH number with $k$ digits and $M=2$. Then $s_2(N)\le k$ and $N\ge 2^{k-1}$. This implies
\begin{equation}
2^{k-1}\le 6k.
\end{equation}
One shows that $6k<2^{k-1}$ for $k\ge 5$ and gets a contradiction.
\end{proof}

\section{Proof of Theorem \ref{thm:2-strong}}\label{sec:4}

\begin{proof} It follows from formula \eqref{eq:9} in the proof of Theorem \ref{thm:2} that:
\begin{equation}\label{eq:1-strong}
b^{k-1}\le (b^2-1)kM.
\end{equation}

We show by induction on the variable $k$ that:
\begin{equation}\label{eq:2-strong}
b^{k-1}>(b^2-1)kM\text{ if }M\ge b^6, k\ge 2\lfloor \log_b M \rfloor+1, b\ge 10
\end{equation}
which together with \eqref{eq:1-strong} ends the proof of Theorem \ref{thm:2-strong} for base $b\ge 10$.

First we show by induction on the variable $M$ that:
\begin{equation}\label{eq:3-strong}
M  >2 b^2(b^2-1)\log_b M+b^2(b^2-1) \text{ if }M\ge b^5, b\ge 10.
\end{equation}

If $M=b^6$ \eqref{eq:3-strong} is equivalent to
\begin{equation}\label{equ:crucial}
b^3+13b(1-b^2)> 0,
\end{equation}
which is true if $b\ge 10$.

Now assume that \eqref{eq:3-strong} is true for a fixed $M$. One has
\begin{equation*}
M+1> 2b^2(b^2-1)\log_b M+b^(b^2-1)+1.
\end{equation*}

To finish the proof of \eqref{eq:3-strong} we still need to check that:
\begin{equation*}
2b^2(b^2-1)\log_b M+b^2(b^2-1)+1\ge 2b^2(b^2-1)\log_b (M+1)+b^2(b^2-1),
\end{equation*}
which after simplifications becomes
\begin{equation*}
1\ge 2b^2(b^2-1)\left ( \log_b(M+1)-\log_b M\right ),
\end{equation*}
which is true due to $M\ge b^5$ and the Mean Value Theorem.

We start the proof of \eqref{eq:2-strong}. In the first step $k=2\lfloor\log_b M\rfloor+1$ and \eqref{eq:2-strong} becomes
\begin{equation}\label{eq:4-strong}
b^{2\lfloor \log_b M\rfloor}>(b^2-1)M(2\lfloor\log_b M\rfloor+1).
\end{equation}

Due to $\log_b M-1\le \lfloor \log_b M\rfloor \le \log_b M$ one has
\begin{equation}\label{eq:5-strong}
\begin{gathered}
b^{2\lfloor\log_b M\rfloor}\ge b^{2(\log_b M-1)}\\
(b^2-1)M(2\lfloor\log_b M\rfloor+1)\le (b^2-1)M(2\log_b M+1).
\end{gathered}
\end{equation}

In order to prove \eqref{eq:4-strong} it is enough to show that
\begin{equation*}
b^{2(\log_b M-1)}> (b^2-1)M(2\log_b M+1),
\end{equation*}
which is equivalent to \eqref{eq:3-strong}. This ends the proof of the first induction step.

Now assume that \eqref{eq:2-strong} is true for fixed $k$ and show that it is true for $k+1$. Due to the induction hypothesis one has that:
\begin{equation*}
b^k\ge b\cdot (b^2-1)kM.
\end{equation*}

To finish the proof of \eqref{eq:2-strong} we still need to check that
\begin{equation*}
b\cdot (b^2-1)kM > b(b^2-1)(k+1)M,
\end{equation*}
which is obviously true.

The proofs of the other cases are similar. The only significant difference appears in \eqref{equ:crucial}. If $3\le b\le 9$, \eqref{equ:crucial} becomes $b^4-15(b^2-1)\ge 0$, which is true. If $b=2$ \eqref{equ:crucial} becomes $b^6>17(b^2-1)$, which is true.
\end{proof}

\section{Proof of Theorem \ref{thm:3}}\label{sec:5}

\begin{proof} As $N$ has $k$ digits one has that:
\begin{equation}\label{eq:15}
N\ge b^{k-1}.
\end{equation}

The largest possible value for $s_b(N)$ is $(b-1)k$. Reversing the order of the digits in an integer increases its value by at most $b$ times. One has that:
\begin{equation}\label{eq:16}
Ms_b(N)\cdot (Ms_b(N))^R\le b(b-1)^2k^2M^2.
\end{equation}

Combining equations \eqref{eq:2}, \eqref{eq:15}, \eqref{eq:16} one has that:
\begin{equation}\label{eq:17}
b^{k-1}\le b(b-1)^2k^2M^2.
\end{equation}

Now we prove by induction on the variable $k$ that:
\begin{equation}\label{eq:18}
b^{k-1} > b(b-1)^2k^2M^2, \text{ for }k\ge M+5, M\ge 1, b\ge 6
\end{equation}
which combined with \eqref{eq:17} ends the proof of Theorem \ref{thm:3} for $b\ge 6$.

In the initial induction step $k=M+5$. The statement in \eqref{eq:18} becomes
\begin{equation}\label{eq:19}
b^{M+4} > b(b-1)^2(M+5)^2M^2, \text{ for }M\ge 1, b\ge 6.
\end{equation}

We prove \eqref{eq:19} by induction on the variable $M$. If $M=1$ \eqref{eq:19} becomes $b^5 > 36b(b-1)^2$, which is true if $b\ge 6$.

Now we assume that \eqref{eq:19} is true for $M$ and prove it for $M+1$. From the induction hypothesis one has that:
\begin{equation}\label{eq:20}
b^{M+5}=b\cdot b^{M+4} > b\cdot b(b-1)^2(M+5)^2M^2.
\end{equation}

In order to finish the proof, we still need to check that:
\begin{equation}\label{eq:21}
b\cdot b(b-1)^2(M+5)^2M^2\ge b(b-1)^2(M+6)^2(M+1)^2
\end{equation}
for $M\ge 1$.

After simplifications, \eqref{eq:21} becomes
\begin{equation}\label{eq:22}
(b-1)M^4+(10b-14)M^3+(25b-61)M^2-84M-36\ge 0,
\end{equation}
which is true for $M\ge 1$ and $b\ge 6$.

This finishes the proof of \eqref{eq:19}.

We continue with the general step in the proof of \eqref{eq:18}. By induction
\begin{equation*}
b^k=b\cdot b^{k-1}>b\cdot b(b-1)^2k^2M^2.
\end{equation*}

To finish the proof of \eqref{eq:18} we still need to check that
\begin{equation*}
b\cdot b(b-1)^2k^2M^2\ge b(b-1)^2(k+1)^2M^2,
\end{equation*}
which after simplifications becomes $(b-1)k^2-2k-1\ge 0$. This is true if $k\ge 1$ and $b\ge 6$.

This finishes the proof of Theorem \ref{thm:3} for $b\ge 6$.

The proof of the case $b=5$ is similar. The only significant changes appear in \eqref{eq:19} and in \eqref{eq:22}. Equation \eqref{eq:19} simplifies to $b^5 > 49(b-1)^2$, which is true for $b= 5$, and \eqref{eq:22}  becomes
\begin{equation*}
(b-1)M^4+(12b-16)M^3+(36b-78)M^2-112M-49\ge 0,
\end{equation*}
which is true if $b=5$.

If $2\le b\le 4$, using $M\ge 1$, the statement in the theorem is true if $k\le 8$. We can assume $k\ge 9$.

If $b=4$ \eqref{eq:17} is still true and gives $4^{k-2}\le 9k^2M^2$. It is easy to show by induction that for $k\ge 9$ one has $4^{k-2}>9k^2(k-8)^2$. If $k\ge M+8$ this implies $4^{k-2}> 9k^2M^2$, a contradiction.

If $b=3$, \eqref{eq:17} is still true and gives $3^{k-2}\le 4k^2M^2$. It is easy to show by induction that for $k\ge 9$ one has $3^{k-2}>4k^2(k-8)^2$. If $k\ge M+8$ this implies $3^{k-2}> 8k^2M^2$, a contradiction.

If $b=2$, \eqref{eq:17} is still true and gives $2^{k-2}\le k^2M^2$. It is easy to show by induction that for $k\ge 9$ one has $2^{k-2}>k^2(k-8)^2$. If $k\ge M+8$ this implies $2^{k-2}> 8k^2M^2$, a contradiction.
\end{proof}

\section{Proof of Theorem \ref{thm:3-new}}\label{sec:6}

\begin{proof} It follows from formula \eqref{eq:17} in the proof of Theorem \ref{thm:3} that:
\begin{equation}\label{eq:1*}
b^{k-1}\le b (b-1)^2k^2M^2.
\end{equation}

We prove by induction on the variable $k$ that:
\begin{equation}\label{eq:2*}
b^{k-1}> b(b-1)^2k^2M^2\text{ for }M\ge b^9, k\ge 3\lfloor \log_b M\rfloor +1, b\ge 9,
\end{equation}
which combined with \eqref{eq:1*} finishes the proof of Theorem \ref{thm:3-new}.

We start showing by induction on $M$ that:
\begin{equation}\label{eq:3*}
M> (b-1)^2b^4(3\log_b M+1)^2\text{for }M\ge b^9, b\ge 9
.
\end{equation}

If $M=b^9$ \eqref{eq:3*} becomes, after cancellations,
\begin{equation}\label{eq:changes13}
b^5>28^2(b-1)^2
\end{equation}
 which is true for $b\ge 9$.

We assume now that \eqref{eq:3*} is true for fixed $M$. We show that it is true for $M+1$.
From the induction hypothesis one has that:
\begin{equation*}
M+1> (b-1)^2b^4(3\log_b M+1)^2+1.
\end{equation*}

To finish the proof of \eqref{eq:3*}, one still needs to check that:
\begin{equation*}
(b-1)^2b^4(3\log_b M+1)^2+1\ge (b-1)^2b^4\left ( 3\log_b(M+1)+1\right )^2,
\end{equation*}
which after algebraic manipulations becomes
\begin{equation}\label{eq:4*}
1\ge b^4(b-1)^2(3\log_b(M+1)-3\log_b M)(3\log_b (M+1)M+2).
\end{equation}

Due to the Mean Value Theorem, \eqref{eq:4*} follows if we show that:
\begin{equation}\label{eq:4*bis}
1\ge 3 b^4(b-1)^2\cdot\frac{1}{M}\left ( 3\log_b(M^2+M)+2\right ).
\end{equation}

Consider the function $g(M)=\frac{1}{M}[3\log_b(M^2+M)+2]$, with derivative:
\begin{equation*}
g'(M)=\frac{\frac{1}{\ln b}\cdot\frac{3}{M^2+M}\cdot(2M+1)M-\left ( 3\log_b(M^2+M)+2\right )}{M^2}.
\end{equation*}

For $M\ge b^9$ the first term in the numerator of $g'(M)$ is $\le 6$ and the absolute value of the second term is $\ge 30$. We conclude that $g'(M)$ is negative and $g(M)$ is decreasing on the interval $[b^9,+\infty)$. The value of
\begin{equation}\label{eq:lt1}
3\cdot b^4(b-1)^2\cdot g(M)
\end{equation}
for $M=b^9$ is larger than $\frac{168(b-1)^2}{b^5}$, which shows that \eqref{eq:4*bis} is true if $b\ge 9$. Consequently \eqref{eq:4*} and \eqref{eq:3*} are true.

We start the proof of \eqref{eq:2*}. In the first step $k=3\lfloor \log_b M\rfloor+1$. Equation \eqref{eq:2*} becomes
\begin{equation}\label{eq:11*}
b^{3\lfloor \log_b M\rfloor}>(b-1)^2b(3\lfloor \log_b M\rfloor+1)^2M^2.
\end{equation}

Due to $\log_b M-1\le \lfloor \log_b M\rfloor \le \log_b M$, \eqref{eq:11*} follows if we prove that
\begin{equation}\label{eq:12*}
b^{3(\log_b M-1)}\ge (b-1)^2bM^2(3\log_b M+1)^2\text{ for }M\ge b^9.
\end{equation}

After algebraic manipulations \eqref{eq:12*} is exactly \eqref{eq:3*}, so it is true. Now we show the general induction step for \eqref{eq:2*}.
Assume \eqref{eq:2*} valid for fixed $M$. Then one has that:
\begin{equation*}
b^k\ge b\cdot b^{k-1}\ge b\cdot (b-1)^2bk^2M^2.
\end{equation*}

To finish we still need to check that:
\begin{equation*}
b\cdot (b-1)^2bk^2M^2\ge (b-1)^2b(k+1)^2M^2,
\end{equation*}
which after simplifications becomes $bk^2\ge k^2+2k+1$ which is true for $k\ge 1$.
This ends the proof of the case $b\ge 9$.

The proofs of the other cases are similar and the only significant changes appear in checking the equation \eqref{eq:changes13} and checking that expression \eqref{eq:lt1} is less than $1$. Due to our assumptions the equation remains valid and the expression is less than $1$.
\end{proof}

\section{Examples of ARH numbers}\label{sec:7}

\begin{table}[h]
\centering
\begin{tabular}{ |c | c |}
\hline
 $M$ & $N$ \\
\hline
 1 & 18, 99 \\
 2 & 12, 33, 66, 99\\
 3 & 99\\
 4 & 99\\
 5 &  11,22,33,44,55,66,77,88,99\\
 6 & \\
 7 & 747\\
\hline
\end{tabular}
\caption{ARH numbers with multipliers 1, 2, 3, 4, 5, 7 and without zero digits.}
\label{t:1}
\end{table}

We list in Table \ref{t:1} small additive multipliers $M$ and the corresponding ARH numbers $N$ without zero digits. Theorem \ref{thm:2} shows that an ARH number with multiplier $6$ has at most $8$ digits. A computer search through all integers with at most 8 digits and all digits different from zero, shows that $6$ is not an additive multiplier for numbers with all digits different from zero. If we allow for zero digits one finds that $909$ is an ARH number with multiplier $6$. A computer search through all integers with at most $11$ digits shows that $9$ is not an additive multiplier. These observations motivate Question \ref{q:6}.

We observe that certain ARH numbers, for example $99$, have several additive multipliers, respectively $1,2,3,4,5$. We also observe that certain multipliers, for example $5$, have associated several ARH numbers, respectively $11,22,33,44,55,66,77,88,99$. The last observation motivates the following definition and questions.

\begin{definition} If $M$ is an additive multiplier in a base $b$, let the \emph{multiplicity} of $M$ be the cardinality of the corresponding set of $b$-ARH numbers.
\end{definition}

\begin{question}\label{q:9} If we fix the multiplicity and the base, is the set of additive multipliers infinite?
\end{question}

\begin{question}\label{q:10-bis} If we fix the base, is the multiplicity of additive multipliers bounded?
\end{question}

\section{Examples of MRH numbers}\label{sec:8}

We list in Table \ref{t:2} small multiplicative multipliers $M$ and the corresponding MRH numbers $N$. Theorem \ref{thm:3} shows that a MRH number with multiplier $3$ has at most $7$ digits. A computer search through all integers with at most 7 digits shows that $3$ is not a multiplicative multiplier. This motivates Question \ref{q:7}.

\begin{table}[h]
\centering
\begin{tabular}{ |c | c | }
\hline
 $M$ & $N$ \\
\hline
 1 & 1, 18, 1458, 1729 \\
 2 & 2268, 736\\
 3 & \\
 4 & 1944, 7744\\
 5 & 71685\\
\hline
\end{tabular}
\caption{MRH numbers with multipliers $1,2,3,4,5$ and without zero digits.}
\label{t:2}
\end{table}

One can also arrange the data as in Table \ref{t:3}, where, for small values of $k$, we list multiplicative multipliers $M$ and the corresponding MRH numbers $N$ with $k$ digits.

\begin{table}
\centering
\begin{tabular}{| c | c | c |  c | c | c |}
\hline
 $k$ & $M$ & $N$ & $k$ & $M$ & $N$\\
\hline
 1 & 1 & 1  & 7 & 22 & 9379678 \\\cmidrule{1-3}
 2 & 1 & 81 &  & 28 &  6527836\\\cmidrule{1-3}
 3 & 2 & 736& & 29 & 9253987 \\\cmidrule{1-3}
 4 & 1 & 1458, 1729 & & 32 & 2892672 \\
   & 2 & 2268 & & 33 &  8673885\\
	 & 4 & 1944, 7744 & &  34 &  7526716 \\\cmidrule{1-3}
 5 & 5 & 71685 & & 38  & 3773932, 6362226 \\
   & 7 & 23632 & & 39  &  5673564 \\
	 & 8 & 94528 & &  41 &  2187391 \\
	 & 9 & 42282 & &  49 &  4274613, 8239644 \\
	 & 14 & 51142 & & 63 & 1821771\\
	 & 23 & 78246 & & 72 & 7651584\\\cmidrule{1-3}
 6 & 12 & 132192  & & 73 & 2895472\\
   &  14 & 188356, 247324 & &  82 & 7651584 \\
   &  19 & 161595  & &  84 &  3252312  \\\cmidrule{4-6}
	 &  21 & 433755, 496692 & 8 &  37 & 13184839\\
	 &  22 & 234256 & &  46 & 11361448 \\
	 &  23 &  685584 & &  48 &  14292288 \\
	 &  26 &  258778 & &  53 & 15437628\\
	 &  27 &  332424 &  &  61 & 15178752\\
	 &  29 &  679354 &  & 66 &  15995232\\
   & 31  & 122512 & &  89 &  7331464\\
	 & 33  &  176418 & & 66 &  15995232\\
	 & 34 &  132192, 751842& &  68 & 11715516\\
	 &  36 & 271188  &&  71 & 16746912 \\
	 &  37 & 215821 &  & 74 & 12419568, 15478432 \\
	 &  38 &  332424 &  & 75 & 19348875 \\
	 & 39 &  145314 & & 76 & 17433792 \\
	 &  44 & 235224 & & 77 & 19552995 \\
      & & & & 78 & 12661272, 22694256 \\
     &&&& 79 & 11437225 \\
	 &&&& 86 & 21371688\\
	 &&&& 89 & 12918439 \\
	 \hline
\end{tabular}
\caption{MRH numbers with $1,2,3,4,5,6, 7, 8$ digits and no zero digits.}
\label{t:3}
\end{table}

We observe from Table \ref{t:3} that certain MRH numbers, for example, $332424$, $132192$, and $3252312$, have several multipliers (respectively $\{27, 38\}$, $\{12, 34\}$, $\{72, 82\})$. We also observe from Table \ref{t:2} that certain multipliers, for example $4$, have associated several MRH numbers, respectively $1944, 7744$. The last observation motivates the following definition and questions.

\begin{definition} If $M$ is a multiplicative multiplier in base $b$, let the \emph{multiplicity} of $M$ be the cardinality of the corresponding set of  $b$-MRH numbers.
\end{definition}

\begin{question}\label{q:12} If we fix the multiplicity and the base, is the set of multiplicative multipliers infinite?
\end{question}

\begin{question}\label{q:13-bis} If we fix the base, is the multiplicity of multiplicative multipliers bounded?
\end{question}

\section{Conclusion}\label{sec:9}

In this paper for any numeration base $b$ we introduce two new classes of integers, $b$-ARH numbers and $b$-MRH numbers. They have properties that generalize a property of the taxicab number $1729$. The second class is a subclass of the class of $b$-Niven numbers. We ask several natural questions about these classes and partially answer some of them. In particular, we show that the class of $b$-ARH numbers is infinite if $b$ is even and that the class of $b$-MRH numbers is infinite if $b$ is odd.

Among the questions left open, the most intriguing is if the set of MRH numbers with all digits different from zero is infinite. One way to attack it is to find an infinity of integers $N$ such that $N=N^R$, $N$ is divisible by $s(N^2)$, and $N^2$ has no digit equal to zero. Then the squares are an infinity of MRH numbers with nonzero digits. Our data shows some examples of such integers.
\begin{itemize}
\item $N^2=188356=434^2, s(N^2)=31| 434$,
\item $N^2=234256=484^2, s(N^2)=22| 484$,
\item $N^2=685584=828^2, s(N^2)=36| 828$.
\end{itemize}

\section{Acknowledgments} The author would like to thank the editor and the referee for valuable comments that helped him write a better paper. He also thanks the OEIS Wiki community for help with posting the sequences A305130 and A305131 on OEIS.

\bigskip
\hrule
\bigskip

\noindent 2010 {\it Mathematics Subject Classification}:
Primary 11B83; Secondary 11B99.

\noindent \emph{Keywords:} base, $b$-Niven number, reversal, additive $b$-Ramanujan-Hardy number, multiplicative $b$-Ramanujan-Hardy number, high degree $b$-Niven number.

\bigskip
\hrule
\bigskip

\noindent (Concerned with sequences
\seqnum{A005349},
 \seqnum{A067030},
\seqnum{A305130}, and
\seqnum{A305131}.)

\bigskip
\hrule
\bigskip

\vspace*{+.1in}
\noindent
Received 2018;
revised version received  .
Published in {\it Journal of Integer Sequences},

\bigskip
\hrule
\bigskip

\noindent
Return to
\htmladdnormallink{Journal of Integer Sequences home page}{http://www.cs.uwaterloo.ca/journals/JIS/}.
\vskip .1in


\begin{thebibliography}{99}

\bibitem{B1} S. Boscaro, Nivenmorphic integers, {\it J. Rec. Math.}, {\bf 28} (1996--1997), 201--205.

\bibitem{B2} E. Bloem, Harshad numbers, {\it J. Rec. Math.}, {\bf 34} (2005), 128.

\bibitem{C} T. Cai, On $2$-Niven numbers and $3$-Niven numbers, {\it Fibonacci Quart.}, {\bf 34} (1996), 118--120.

\bibitem{CK} C. N. Cooper and R. E. Kennedy, On consecutive Niven numbers, {\it Fibonacci Quart.}, {\bf 21} (1993), 146--151.

\bibitem{D} W. L. McDaniel, The existence of infinitely many $k$-Smith numbers, {\it Fibonacci Quart.}, {\bf 25} (1987) 76--80.

\bibitem{KD} J. M. De Koninck and N. Doyon, Large and small gaps between consecutive Niven numbers,
{\it J. Integer Seq.}, {\bf 6} (2003), Article 03.2.5.

\bibitem{G} H. G. Grundman, Sequences of consecutive Niven numbers, {\it Fibonacci Quart.}, {\bf 32} (1994), 174--175.

\bibitem{FILS} H. Fredricksen, E. J. Iona\c scu, F. Luca, and P. St\u anic\u a, Remarks on a sequence of minimal Niven numbers, In: {\it Sequences, Subsequences, and Consequences},  Lec. Notes in Comp. Sci., Vol 4893.
Springer, 2007, pp. \ 162--168.

 \bibitem{Hardy} G. H. Hardy, {\it Ramanujan: Three Lectures on Subjects Suggested by his Life and Work}, Chelsea, 1999.

\bibitem{online} N. J. A. Sloane, {\it The On-Line Encyclopedia of Integer Sequences}, \url{http://oeis.org}.

\bibitem{N} V. Ni\c tic\u a, High degree $b$-Niven numbers,  preprint, 2018, \url{http://arxiv.org/abs/1807.02573}.

\end{thebibliography}
\end{document}